\documentclass[a4paper,12pt]{amsart}
\usepackage{amssymb}
\usepackage{latexsym}
\usepackage{amsmath}
\usepackage{enumerate}
\usepackage{amsmath, hyperref}
\usepackage{url} 
\usepackage{appendix}
\usepackage{tikz}
\usepackage{geometry}
\usepackage[all,pdf]{xy}
\usepackage{enumitem}
\setlist{nosep}

\theoremstyle{definition}

\begin{document}

\title{On the depth of Wiles's  proof  of Fermat's Last Theorem}

\begin{abstract}
What constitutes depth in a mathematical proof? This paper addresses this foundational question through a close case study of Andrew Wiles's proof of Fermat's Last Theorem. We propose a five-criteria framework—Difficulty, Originality, Fruitfulness, Unity, and Explanatory Power—to analyze proof depth. By critically examining the limitations of each criterion and the intricate interrelations among them, we show that they do not merely coexist but capture genuinely complementary aspects of proof depth. 

Applying this framework to Wiles's proof, we then demonstrate that its celebrated depth does not reside in any single virtue, but rather emerges from a synthesis of all five criteria. Our analysis thus offers not only a novel lens for understanding Wiles's achievement, but also a flexible framework for evaluating depth across mathematical proofs more broadly. 
\end{abstract}

%\subjclass[2010]{03F40, 03F30, 03F25}


\author{Yong Cheng and Colin McLarty}
\date{}

\maketitle

\section{Introduction}\label{sec1}

The term ``mathematical depth" is highly valued in the mathematical community as a measure of a mathematician's work, yet it remains an elusive concept. Since 2003, the Abel Prize laureates' work has been called ``deep" nine times [Holden and Piene, 2014]; a 2014 special issue of \emph{Philosophia Mathematica} was dedicated to the topic [Ernst et al., 2014].\footnote{The papers in that issue vary in methodology, focus, and conclusions, though they are largely complementary.} 

Numerous criteria for assessing mathematical depth have been proposed, primarily focusing on the depth of mathematical theorems [Ernst et al., 2014; Cheng, 2022]. Currently, there are no universally accepted criteria for characterizing mathematical depth. Existing accounts have emphasized different features, such as proof difficulty or complexity, theorem fruitfulness, or the capacity to unify disparate domains.\footnote{See Urquhart [2014], Stillwell [2014], Arana [2015], Tappenden [2012], and Ernst et al. [2014].} The literature reflects a strong consensus that there is no single way to achieve depth [Ernst et al., 2014]; rather, mathematical depth appears to be a multifaceted concept with multiple, irreducible manifestations—a perspective we develop below under the rubric of ``family resemblance" (following Wittgenstein [2009]).

The special issue in \emph{Philosophia Mathematica} concludes that mathematical depth is a complex, multifaceted concept that resists simple definition, yet it remains a rich and legitimate subject for philosophical and methodological inquiry, warranting further study [Ernst et al., 2014].

The concept of mathematical depth can apply to various entities, including theorems, proofs, definitions, axioms, fields, methods, and conjectures. We distinguish between the depth of a mathematical proof and that of a theorem in two respects. First, a theorem's being deep does not guarantee that every proof of it is deep (a theorem may admit multiple proofs that differ significantly in nature, raising the question of which, if any, are deep or which is the deepest). Second, the depth of a proof does not automatically transfer to the theorem it proves.

This paper focuses on the depth of Wiles's proof of Fermat's Last Theorem (FLT), not on the theorem itself.\footnote{Throughout the paper, we use FLT to refer to Fermat's Last Theorem.} Three reasons motivate this focus: (i) no fundamentally different proof of FLT exists; (ii) FLT itself, unlike the Taniyama–Shimura–Weil conjecture, yields few consequences [McLarty, 2024]; and (iii) the depth of proofs has received less attention than the depth of theorems.

This paper does not aim to define proof depth, but to establish a synthetic framework for its analysis. We propose a five-criteria framework for analyzing the depth of mathematical proofs:\emph{Difficulty}, \emph{Originality}, \emph{Fruitfulness}, \emph{Unity}, and \emph{Explanatory Power}. In Sections \ref{sec 3.1} to \ref{sec 3.5}, we examine each criterion, discuss potential drawbacks, and explore their interrelationships.  Our framework synthesizes these five criteria, capturing proof depth through their distinct roles and mutual dependencies.

We apply our five-criteria framework to Wiles's proof of FLT. Wiles's proof has been widely described as deep—e.g., as ``stunning" by the Norwegian Academy of Science and Letters in awarding the 2016 Abel Prize to Wiles, as ``one of the deepest achievements in the history of mathematics" by the American Mathematical Society,\footnote{American Mathematical Society, ``Fermat's Last Theorem (2-Volume Set)", AMS Bookstore, 2014, available at \url{https://bookstore.ams.org/mmono-243-245}.} and as ``extremely deep" by number theorist Andrew Granville.\footnote{A. Granville, ``Fermat's Last Theorem announcement", Fermat News, Newton Institute, Cambridge, 1993, available at \url{http://felix.unife.it/Root/d-Mathematics/d-Number-theory/d-News/t-Fermat-R.I.P}.} We examine the philosophical basis for this judgment.

We explore the question: What constitutes the depth of Wiles's proof? We argue that its depth is best understood as a synthesis of the five proposed criteria. By `synthesis' we mean not merely that Wiles's proof scores highly on each criterion taken separately, but that in this proof the criteria interact and reinforce each other in an integrated way. This mutual reinforcement is not a mere correlation; it consists of concrete dependencies among the criteria—for example, unity enabling explanatory power, originality and unity generating fruitfulness. We develop the family-resemblance thesis in Section \ref{literature review}, where we explicitly adopt it as the philosophical grounding for our five-criteria framework.

The paper proceeds as follows. Section \ref{literature review} reviews the literature. Section \ref{sec2} outlines Wiles's proof of FLT. Section \ref{sec3} presents our five-criteria framework. Section \ref{sec4} applies this framework to assess the depth of Wiles's proof.  Section \ref{sec5} concludes.

\section{Review of Current Literature}\label{literature review}

We review existing accounts of mathematical depth, pinpoint their limitations, and thereby motivate our five-criteria framework.

Arana [2015] distinguishes four views of mathematical depth: the Genetic View (depth tied to the mathematician's talent), the Evidentialist View (depth tied to proof properties), the Consequentialist View (depth tied to fruitful consequences), and the Cosmological View (depth tied to revealing structure).\footnote{See Arana [2015] for detailed exposition. The Genetic View posits that a theorem is deep if proven by exceptionally talented mathematicians—a view criticized for its subjectivity and its potential to classify trivial results as deep ([Arana, 2015]). The Evidentialist View claims that a theorem's depth is tied to properties of its proof(s) ([Arana, 2015]). The Consequentialist View asserts that a theorem is deep if it has profound, fruitful consequences ([Arana, 2015]). The Cosmological View claims that a theorem is deep if it reveals significant order or structure ([Arana, 2015]). Arana critiques these views using Szemer\'{e}di's Theorem as a case study, arguing that all of them have problems as measures of depth and that depth is a complex, multifaceted concept that resists simple definition.}  His analysis focuses on the depth of a theorem rather than that of a proof. Arana's paper [2015] is primarily critical, identifying problems with existing views, concluding without offering a synthesized perspective or a clear path forward. 
 
Stillwell [2014] proposes three types of depth—Historical, Foundational, and Logical—each emphasizing different aspects of mathematical practice.\footnote{Historical Depth refers to the cumulative effort of mathematicians over generations, indicating that the underlying ideas were not easily accessible ([Stillwell, 2014]). Foundational Depth signifies that a theorem is fundamental and fruitful, supporting a large body of further mathematics ([Stillwell, 2014]). Logical Depth is a more formal measure, often related to the length, difficulty, and laboriousness of the shortest or only known proof ([Stillwell, 2014; Arana, 2015; Urquhart, 2014]).} However, Stillwell's work, while insightful, remains largely taxonomic and lacks a synthetic framework to show how the three types of depth interact or how they might be applied comparatively. 
 
Gray [2014] characterizes deep mathematics as hidden/difficult,\footnote{That is, not immediately accessible; its truths are not surface-level ([Gray, 2014]).} structural,\footnote{It reveals underlying coherence and necessitates a reorganization of existing knowledge ([Gray, 2014]).} and explanatory.\footnote{It does not merely prove a fact; it explains why it is true ([Gray, 2014]).} However, Gray's paper is firmly historical and offers limited discussion of ``depth" in a modern context. While it establishes a powerful paradigm, it does not explore how well the Gaussian definition of depth holds up in 20th and 21st-century mathematics.

The afterword in Ernst et al. [2014] identifies five promising criteria for capturing mathematical depth, referred to as ``Candidate Criteria": (1) ties together apparently disparate fields; (2) involves impurity;\footnote{That is, using concepts from outside the problem's native domain or from ``higher" conceptual strata. An ``impure" proof often signifies depth by revealing unexpected connections ([Arana, 2015; Urquhart, 2014]).} (3) finds order in chaos;\footnote{That is,  reveals hidden structure in apparent disorder.} (4) exhibits organizational or explanatory power; (5) transforms a field or opens a new one.\footnote{No examples of deep mathematical theorems that fail to meet any of these criteria have been found ([Ernst et al., 2014]).} 
These criteria summarize discussions of mathematical depth in various papers from the special issue of \emph{Philosophia Mathematica}. However, no single paper has analyzed mathematical depth using all five criteria.

None of the existing accounts provides a set of criteria that are simultaneously (i) applicable to proofs, (ii) mutually interacting rather than isolated, and (iii) accompanied by a clear heuristic for application. These limitations motivate our five-criteria framework for proof depth (Section \ref{sec3}), which synthesizes existing insights: Difficulty draws on Evidentialist, Logical, and Historical views; Originality refines the Genetic view; Fruitfulness aligns with Consequentialist and Foundational views; Unity extends the Cosmological view; and Explanatory Power builds on Lange [2014], Gray [2014], and Steiner [1978].

Ernst et al. [2014] conclude that there is no single kind of mathematical depth; rather, there are multiple, irreducible ways of being deep, and these different ways often come apart—a view that aligns with Wittgenstein's later notion of family resemblance ([Wittgenstein, 2009]). Just as games share overlapping similarities rather than a single common property, mathematical depth is constituted by a cluster of features that may vary across cases.\footnote{See Wittgenstein [2009] on ``game" as a family-resemblance concept.} The five criteria we propose are those that cluster together in paradigmatic cases—they are not an invariant set required for every deep proof. 

We explicitly adopt this family-resemblance perspective. It implies that any attempt to define ``mathematical depth" by a set of necessary and sufficient conditions is misguided. Our five-criteria framework is therefore not a definition but a heuristic tool for systematic comparison and philosophical analysis. The five criteria are features that cluster together in paradigmatic cases of deep proofs—such as Wiles' proof of FLT—and they serve as a structured vocabulary for articulating why a proof is considered deep. As we argue in Section \ref{sec 3.7}, none of the five criteria is strictly necessary, nor is any one of them sufficient on its own. Their value lies not in providing a checklist for depth, but in capturing the multifaceted nature of the concept and facilitating systematic comparison. 

In the course of our argument we also answer a potential objection: if depth is a family-resemblance concept, why propose a fixed set of five criteria?  The answer is that the criteria are not intended to be exhaustive or final. They are salient features extracted from the literature and from mathematical practice, and they have proven their utility in the paradigmatic case of Wiles' proof as we will show in Section \ref{sec4}. Other proofs may bring other features to the fore, and the framework is open to revision. What matters is that the criteria provide a common ground for comparing different instances of depth, not that they capture every possible manifestation.

\section{An Outline of Wiles's Proof of Fermat's Last Theorem}\label{sec2}

We provide an overview of Wiles's proof of FLT.\footnote{A vast literature exists on Wiles's proof, catering to various levels of sophistication. Notably, McLarty [2024] targets philosophers of mathematics, while Mazur [1991] offers insights for mathematicians outside number theory. Darmon et al. [1997] is also highly regarded.} FLT states that for any integer $n> 2$, there are no positive integers $a, b$ and $c$ such that $a^n + b^n = c^n$.  First conjectured by Fermat in 1637, it was proved by Andrew Wiles in 1995.\footnote{Pierre de Fermat first proposed this conjecture in the margins of his book \emph{Arithmetica} around 1637.}
 
We begin by introducing key concepts. \emph{Elliptic curves} are defined by cubic equations in two variables (e.g., $y^2 = x^3 + ax + b$) and are fundamental objects in number theory with rich geometric and algebraic structures. \emph{Modular forms} are highly symmetric complex functions whose Fourier coefficients encode arithmetic data and satisfy numerous symmetry conditions. Before Wiles's work, modular forms primarily resided in complex analysis. An elliptic curve is \emph{modular} if its arithmetic properties  match the Fourier coefficients of a modular form. A \emph{Galois representation} captures the intricate symmetries of number systems via linear algebra.\footnote{The Absolute Galois Group of $Q$, denoted $G_{Q}$, consists of all symmetries of algebraic numbers that fix the rational numbers. A Galois representation is a continuous group homomorphism $\rho: G_{Q} \rightarrow GL_n(F)$, where $GL_n(F)$ is the group of invertible $n \times n$ matrices over a field $F$.} This representation serves as a ``translation device" that converts abstract algebraic problems into concrete matrix problems.

Around 1955, Japanese mathematicians G. Shimura and Y. Taniyama conjectured that every elliptic curve defined over the field of rational numbers has an associated modular form. This conjecture gained prominence in the West when number theorist Andr\'{e} Weil provided supporting evidence, although he did not prove it.\footnote{By the 1980s, substantial evidence had accumulated, and many papers explored the implications of this conjecture, though it remained unproven at that time.} The conjecture, now known as the \emph{Taniyama-Shimura-Weil Conjecture} (or \emph{Modularity Theorem}), asserts that every elliptic curve over the rational numbers is modular. This conjecture is fundamentally about arithmetic due to its foundational connections to elliptic curves and modular forms, both central to number theory.\footnote{As Mazur [1991] noted, a subtle issue in formulating the Taniyama-Shimura-Weil Conjecture is that, while it is ``about arithmetic", it can be expressed in various ways. For instance, one interpretation connects it to integral transforms in complex variable theory, while another relates it to geometry ([Arana, 2025, p. 36]). Its diverse expressions—through integral transforms, geometric interpretations, algebraic geometry, and representation theory—highlight its richness and the depth of its implications across different mathematical fields.}

Wiles did not prove FLT directly; instead, he proved the Taniyama-Shimura-Weil Conjecture for a specific class of elliptic curves. His approach was groundbreaking, linking this conjecture to FLT through the work of Gerhard Frey, further refined by Jean-Pierre Serre and Ken Ribet. To date, no direct method for proving FLT without involving the Taniyama-Shimura-Weil Conjecture has been discovered. The foundational breakthrough came from Gerhard Frey's innovative idea: he proposed that if a counterexample to FLT exists, one could construct a special elliptic curve—now known as the \emph{Frey curve}—from it.\footnote{Hellegouarch [1971] and Demyanenko [1971] connected FLT to the geometry of special elliptic curves. Mazur [1977, 1978], not aiming at FLT, found key facts on modular elliptic curves. Hellegouarch, Demyanenko, and Mazur all made breakthroughs leading to G. Frey's discovery.} Initially, the Taniyama-Shimura-Weil conjecture and FLT appeared unrelated. Frey's insight transformed the arithmetic problem of FLT into a geometric one: whether this specific Frey curve could be modular. Ken Ribet, building on Jean-Pierre Serre's work, provided the necessary rigor to Frey's vision, proving that if a Frey curve were modular, a sophisticated ``level-lowering" argument would lead to a logical contradiction. Thus, Frey curves cannot be modular, establishing the critical link: the Taniyama-Shimura-Weil conjecture implies FLT, as a counterexample to FLT would yield a non-modular elliptic curve, contradicting the conjecture.\footnote{In 1985, J.P. Serre partially proved that Frey curves are not modular [Cornell et al., 1997]. Serre did not provide a complete proof, and the missing part is known as the epsilon conjecture. In 1986, Ribet successfully proved the epsilon conjecture based on Serre's partial proof, demonstrating that if the Galois representation associated with an elliptic curve has certain properties (which the Frey curve possesses), then the curve cannot be modular [Ribet, 1990].}

Wiles proved FLT by contradiction. He first assumed FLT is false, i.e., there exists a non-zero solution to $a^p + b^p = c^p$ for some odd prime $p > 2$. In the mid-1980s, G. Frey had the insight to construct an elliptic curve from this hypothetical solution: $y^2 = x(x - a^p)(x + b^p)$, known as the ``Frey Curve" [Frey, 1986]. This curve is semi-stable and possesses special properties.\footnote{The ``semi-stable" condition Wiles used is a technical property relating to the number of solutions the curve has modulo primes.} However, Frey did not provide a complete proof that the Frey curve is non-modular. Ribet proved that if such a Frey curve exists, it must be non-modular. Thus, from Ribet's theorem, the Frey curve is non-modular.

Next, Wiles proved the Taniyama-Shimura-Weil conjecture for semi-stable elliptic curves, demonstrating that every semi-stable elliptic curve over the rational numbers is modular. This led to a contradiction: if a counterexample to FLT exists, then the Frey curve must be semi-stable non-modular (by Ribet), yet Wiles proved that every semi-stable elliptic curve (including the Frey curve) is modular. Therefore, a counterexample to FLT cannot exist, confirming that FLT is true.

The key to proving FLT lies in establishing that the Taniyama-Shimura-Weil conjecture holds for semi-stable elliptic curves, which is sufficient to encompass the Frey curve. On October 24, 1994, Wiles submitted two manuscripts titled ``Modular Elliptic Curves and Fermat's Last Theorem" and ``Ring Theoretic Properties of Certain Hecke Algebras", the latter co-authored with R. Taylor (see Wiles [1995], Taylor-Wiles [1995]). Both papers were published in the Annals of Mathematics in May 1995, proving that the Taniyama-Shimura-Weil conjecture holds for semi-stable elliptic curves and thereby resolving FLT.

Wiles's strategy to prove the Taniyama-Shimura-Weil conjecture for semi-stable elliptic curves involved establishing the \emph{Modularity Lifting Theorem}, which states that if the Galois representation of a semi-stable elliptic curve modulo some prime $p$ is modular, then the curve itself is modular. This theorem allows for ``lifting" modularity from a simpler representation modulo a prime to the full structure. Wiles proved the Taniyama-Shimura-Weil conjecture for semi-stable elliptic curves using an induction argument on prime numbers: he first established the Modularity Lifting Theorem and then demonstrated that the Galois representation of a semi-stable elliptic curve modulo some prime $p$ is modular.

Wiles's approach to proving the Modularity Lifting Theorem involved establishing the \emph{$R=T$ Theorem}, which asserts that the \emph{deformation ring} $R$ is isomorphic to the \emph{Hecke algebra} $T$. Here,  $R$ serves as a ``parameter space" for all possible Galois representations lifting a base representation, while $T$ serves as a ``parameter space" for modular forms, encoding their arithmetical properties. The isomorphism $R \cong T$ indicates that every Galois representation in $R$ corresponds to a modular form in $T$, thereby establishing modularity. Proving the $R=T$ Theorem directly establishes the Modularity Lifting Theorem, demonstrating that if the base representation is modular, all its ``nice" lifts are also modular.

Proving \(R \cong T\) directly is technically difficult. Wiles's \(R=T\) strategy forces this isomorphism via a ``controlled comparison'': it constructs a surjection \(\phi: R \to T\) (associating modular forms with Galois deformations, though possibly with ``extra'' ones in \(R\)), and then proves injectivity to eliminate the extraneous cases.\footnote{Injectivity is achieved via the Taylor--Wiles method: auxiliary primes \(Q\) make the isomorphism \(R_Q \cong T_Q\) accessible in a smoothed setting, and a limiting argument over \(Q\) recovers the original \(R \cong T\).} 

The Modularity Lifting Theorem follows directly from the $R=T$ Theorem. To prove the Taniyama-Shimura-Weil conjecture for semi-stable elliptic curves, it suffices to show that the Galois representation of a semi-stable elliptic curve modulo some prime $p$  is modular. If the Galois representation of a semi-stable elliptic curve modulo $p=3$ is irreducible, it has been known since around 1980 that its Galois representation is always modular. However, proving that the Galois representation is modular when the representation modulo $p=3$ is reducible poses significant challenges.

To address this difficulty, Wiles employed the so-called ``$3$–$5$ switch" trick. When the representation of a semi-stable elliptic curve $E$  modulo  $p=3$  is reducible, it becomes easier to work with  $p=5$  and use the Modularity Lifting Theorem to demonstrate that the representation of  $E$  modulo  $p=5$  is modular, rather than proving directly that the reducible representation of  $E$  modulo  $p=3$  is modular.\footnote{Wiles resolves the reducible mod~3 case via a switch to~5. If the mod~5 Galois representation \(E[5]\) is reducible, modularity follows directly. If not, he constructs a semi-stable \(F\) with \(F[3]\) irreducible and \(E[5] \cong F[5]\). By Langlands--Tunnell, \(F[3]\) is modular, hence so is \(F\), and the mod~5 isomorphism transmits this modularity to \(E\).}

Thus, Wiles proved that the Galois representation of a semi-stable elliptic curve modulo some prime $p$ is modular. By applying the Modularity Lifting Theorem, it follows that every semi-stable elliptic curve over the rational numbers is modular, thereby completing the proof of the Taniyama-Shimura-Weil conjecture for semi-stable elliptic curves.

\section{A Five-Criteria Framework for Analyzing Proof Depth}\label{sec3}

We introduce a heuristic five-criteria framework—Difficulty, Originality, Fruitfulness, Unity, and Explanatory Power—for analyzing the depth of mathematical proofs. The framework synthesizes various proposed measures of depth from the literature, but it is not a definition of proof depth; rather, it is a tool for comparative assessment. 

Sections \ref{sec 3.1}–\ref{sec 3.5} examine each criterion and its potential drawbacks. Section \ref{sec 3.6} compares the criteria by highlighting their key differences and interconnections. Section \ref{sec 3.7} evaluates whether the criteria are necessary and sufficient, and Section \ref{sec 3.8} addresses their subjectivity and objectivity. Throughout, our evaluation focuses on the relevant mathematical community rather than on individual authors or readers, thereby minimizing subjectivity rooted in personal interests or abilities.\footnote{For a proof in a given field, the relevant mathematical community comprises the active researchers whose specialized expertise bears directly on its methods and content.}

\subsection{The Criterion of Difficulty}\label{sec 3.1}

Difficulty is the intrinsic challenge a proof presents to an agent with sufficient background, comprising both the conceptual obstacle of uncovering underlying ideas and the technical complexity of the reasoning itself. It is a property of the proof as presented, not merely of its historical discovery.\footnote{Our criterion synthesizes Arana's Evidentialist properties, Urquhart's complexity and non-obviousness, Gray's hiddenness, and Stillwell's Logical metrics (Section \ref{literature review}).}
   
Difficulty has two distinct but complementary aspects. \emph{Conceptual difficulty} concerns the intellectual obstacles involved in discovering the necessary principles or framing shifts—for instance, recognizing that Fermat's Last Theorem could be recast as a statement about elliptic curves and modular forms. It often demands the invention of new ideas (originality is taken up in Section~\ref{sec 3.2}).
      
\emph{Technical difficulty}, by contrast, concerns the complexity of the proof as laid out: overall length, the number and diversity of fields employed, the level of prerequisites required, and structural intricacy (e.g., multi-step, layered reasoning). A proof drawing on fields with divergent methodologies—say, algebraic geometry and analytic number theory—is generally more technically difficult than one confined to a single domain.

These dimensions are analytically distinct but often co-occur. A proof may be conceptually difficult yet technically straightforward (Cantor's diagonal argument: a short construction that reorients our understanding of infinity). Conversely, a proof may be technically massive yet conceptually shallow (the Appel--Haken proof of the Four Color Theorem: an exhaustive case analysis within an existing framework). Wiles's proof of FLT excels in both, exhibiting deep paradigm shifts alongside extraordinary technical intricacy (see Section \ref{sec4}). We take conceptual difficulty to be the more philosophically interesting dimension, as it captures the kind of creative breakthrough that can fundamentally advance—and in paradigmatic cases transform—a field.

Formal measures can serve as objective proxies for technical difficulty, yet they diverge from what mathematicians ordinarily mean by ``difficulty''. Proof complexity, computational complexity, and logical complexity (as studied in reverse mathematics) offer objective metrics, but they are largely independent of perceived difficulty.\footnote{Proof complexity studies the minimal size of a proof within a formal system, offering insights into the proof's intrinsic logical size. Computational complexity analyzes the time and space resources required by algorithms related to the proof. Logical complexity (investigated in Reverse Mathematics) examines the strength of the axioms needed to formalize the proof.} A proof may have low proof complexity and weak axioms yet be conceptually arduous (Cantor), or it may be case-heavy and axiomatically strong while resting on a conceptually flat idea (Appel--Haken). Likewise, the fact that one proof requires stronger set-theoretic axioms than another does not imply that it is more difficult in the sense that matters to mathematical practice.

It is also crucial to distinguish difficulty from the historical process of discovery. Historical evidence—the duration a problem remained open, the number of failed attempts, the cumulative labour of generations—may \emph{suggest} that a proof is difficult, but it is neither necessary nor sufficient for genuine difficulty. A problem may resist solution for centuries simply because the right idea was not yet conceived; conversely, a highly technical proof may be found quickly once the requisite machinery is in place. Thus, historical evidence serves at best as a heuristic indicator. When we assess Wiles's proof in Section~\ref{sec4}, we will cite the 358-year gap and the many failed attempts as supporting evidence for its difficulty, not as a separate type of difficulty.
 
Despite its merits, the difficulty criterion has three main drawbacks: subjectivity across evaluators and over time, relativity across mathematical subfields, and a risk of undervaluing elegant simplicity. These do not invalidate the criterion, but they require assessments to be relativised to community and context, and they demand that simplicity be recognised as an independent virtue (for interactions with other criteria, see Section~\ref{sec 3.6}). 

\subsection{The Criterion of Originality}\label{sec 3.2}

The originality of a proof consists in the introduction of ideas or methods not previously conceived or implemented by the relevant mathematical community at the time of presentation.\footnote{Our criterion refines Urquhart's non-obviousness [2014] and Arana's Genetic View [2015] by focusing on the novelty of ideas and methods introduced in the proof itself.} As the community's knowledge evolves, a proof initially deemed original may be reassessed if similar ideas are discovered to have been present elsewhere; conversely, iterations of known methods that yield no new conception do not count as original.

We distinguish between \emph{invention} and \emph{creative synthesis}. 
Invention introduces genuinely unprecedented ideas or methods. 
Creative synthesis, by contrast, combines existing ideas or methods to solve a problem that neither could address alone; its originality lies not in the novelty of its components but in the unexpected integration that makes the solution possible. Overemphasizing invention risks undervaluing this kind of originality, as a proof may introduce no fundamentally new concept yet achieve its result through a brilliant and unforeseen combination of established tools. Indicators of such originality include novel concepts, unexpected integrations, and significant problem reframings.

The originality criterion has several drawbacks. Determining true originality is historically demanding, as an idea may have been developed independently elsewhere; judgments are relative to a mathematical community, so a proof may be considered novel in one context but derivative in another; and an overemphasis on originality may incentivize contrived, ad-hoc methods when simpler approaches would be more effective. For how this criterion interacts with the others, see Section~\ref{sec 3.6}.   

\subsection{The Criterion of Fruitfulness}\label{sec 3.3}

The fruitfulness of a proof consists in the breadth of its consequences—the range of implications, results, or applications extending beyond its initial context—and the applicability of its methods to other domains.\footnote{This criterion draws on Arana's Consequentialist View [2015],  Stillwell's Foundational Depth [2014], Urquhart's field-impact emphasis [2014], Tappenden's notion of productive richness [2012], and Ernst et al.'s field-transformation criterion [2014].  Tappenden [2012] defines fruitfulness as the capacity of a concept, definition, or theory to generate novel discoveries and connections.  His central contribution is the notion of ``productive richness", although it lacks a precise definition. Our criterion for fruitfulness encompasses Tappenden's concept while offering a broader interpretation.} A fruitful proof catalyzes further exploration, enabling its techniques to be generalized or adapted, thereby generating new theorems, tools, and cross-domain applications. These indicators collectively signal a proof's significance and its capacity to inspire ongoing mathematical inquiry.

The fruitfulness criterion has several drawbacks. Fruitfulness is often recognized only in hindsight, as with Galois's work on group theory, which was largely ignored during his lifetime. Assessments vary across communities: number theorists may find Wiles's proof fruitful, while financial mathematicians may not. There is also a field bias favoring prominent areas like number theory or geometry over equally significant work in niche fields. Moreover, the boundary between breadth of consequences and applicability of methods is vague, as the two are interconnected rather than mutually exclusive. Finally, what counts as a sufficient range of applications for a proof to be deemed fruitful remains undefined. For how this criterion interacts with the others, see Section~\ref{sec 3.6}.
    
\subsection{The Criterion of Unity}\label{sec 3.4}

The unity of a proof consists in the connections it establishes among at least two distinct mathematical domains—objects, theories, or fields—illuminating underlying principles, structures, or methodologies that transcend individual areas. Unity is characterized not by the number of domains involved, but by the theoretical importance and cognitive salience of the connections revealed. A connection is theoretically important if it enables the transfer of problems, methods, or results from one domain to another, or if it reveals previously unrecognized structural isomorphisms. Such connections enhance our understanding of the mathematical landscape without requiring an antecedent judgment of the proof's ``depth".\footnote{One might worry that ``theoretical importance" or ``cognitive salience" merely swaps one vague term for another. Yet these are pre-theoretic anchors identifiable independently of depth—e.g., by cross-field citations, subsequent theorems, or historical testimony. Since unity is only one of five interacting criteria, no circular definition of depth is invoked.}   Indicators of unity include cross-domain principles, structural connections, and methodological syntheses. Our criterion of unity emphasizes a proof's integrative power—its ability to unify disparate mathematical domains—and draws on Arana's cosmological view of depth as revealing order and structure [2015], Ernst et al.'s candidate criteria of cross-field unification and finding order in chaos [2014], and the notion of impurity [Arana, 2015; Ernst et al., 2014].

Proofs that unify distinct domains frequently exhibit impurity—using tools from outside the theorem's ``native" domain. Arana [2025] systematically examines purity, defining a proof as pure if it relies only on elements that are ``close", ``intrinsic", or ``native" to the theorem, avoiding ``extraneous", ``foreign", or ``remote" concepts.\footnote{Arana [2025] argues that purity comprises a cluster of epistemic virtues, identifying five types: geographical, topical, syntactic, logical, and elemental purity.} Despite his emphasis on purity, Arana [2025] acknowledges that impurity can have virtues—unification, discovery, and efficiency among them. Wiles's proof is a paradigmatic case of fruitful impurity: it imported the modularity theorem and Galois representations into number theory, achieving a unification that a pure proof could not.

The unity criterion has several drawbacks. It may reward superficial connections that link domains artificially without genuine insight; the value of a connection lies in deepening understanding, not merely in the act of connecting. What counts as ``distinct" depends on human categorization—objects from two domains may be seen as belonging to a more abstract, higher-order domain. An emphasis on unity may undervalue specialization, favoring horizontal links over vertical breakthroughs within a single domain. It may also devalue purity: a unifying proof can be more abstract and harder to grasp than a direct, pure proof, yet this criterion risks diminishing the achievement of finding a pure proof. Finally, the criterion is qualitative, making it difficult to specify how closely domains must be tied for a proof to be considered unified. For how this criterion interacts with the others, see Section~\ref{sec 3.6}.

\subsection{The criterion of explanatory power}\label{sec 3.5}

The explanatory power of a proof consists in the insights it provides into why the theorem it establishes is true, going beyond mere verification to offer genuine understanding.\footnote{This criterion builds on Lange's emphasis on explanatory power and answering ``why" questions [2014], Gray's explanatory aspect [2014], Steiner's characterizing property [1978], and Ernst et al.'s ``explanatory/organizational power" candidate criterion [2014].} An explanatory proof addresses ``why" questions by connecting results to foundational principles or uncovering essential mechanisms. Steiner [1978] characterizes such proofs as relying on a ``characterizing property"—a property unique to an entity within a relevant domain—and proposes a ``deformation test": a proof is explanatory if, by substituting a different characterizing property for a related object in the same family, we can generate new related theorems.\footnote{Steiner's deformation test operationalizes explanatory power, though it has faced objections (see Hafner and Mancosu [2005]). We adopt it here as a working heuristic because it captures an important aspect of what mathematicians often cite as explanatory in the context of Wiles's proof.} Lange [2014] argues that mathematical depth is linked to a proof's ability to answer ``why" questions, though his account risks circularity: a proof is deep because it is explanatory, and explanatory because it uses a salient feature, deemed salient precisely because it aligns with deep proofs. To avoid this circle, Lange must provide an independent criterion for salience.

We note that although explanatory proofs are sometimes associated with an ``Aha!" moment of insight, this phenomenological marker is neither necessary nor sufficient—it varies across individuals and communities. Our analysis therefore centers on Steiner's characterizing property and on the proof's capacity to answer ``why" questions. Indicators of explanatory power include conceptual reasoning, illumination of essence, characterizing properties, and the revelation of underlying patterns.
   
The explanatory power criterion has several drawbacks. It is subjective: what one mathematician finds illuminating, another may find confusing—a category theorist, a combinatorialist, and an analyst may each prefer a different proof of the same theorem. It also suffers from vagueness: there is no precise definition of what constitutes explanatory insight, which may involve reducing a theorem to intuitive axioms, providing a compelling narrative, offering visual intuition, presenting a constructive proof over a proof by contradiction, or articulating a characterizing property. This lack of clarity complicates consistent application of the criterion. For how this criterion interacts with the others, see Section~\ref{sec 3.6}.
 
In summary, the five criteria—Difficulty, Originality, Fruitfulness, Unity, and Explanatory Power—synthesize diverse perspectives on mathematical depth from the literature, offering a framework for analyzing the depth of mathematical proofs despite their potential drawbacks. 

\subsection{Interrelations Among the Five Criteria}\label{sec 3.6}

The five criteria are distinct but systematically interdependent. Their overlaps, tensions, and causal dependencies—rather than any checklist function—give the framework its analytical purchase. 
     
\begin{itemize}
\item Difficulty correlates with originality and unity, but neither entails the other: a proof may be difficult yet unoriginal (brute-force computation) or disunified (a patchwork of unrelated lemmas), while originality or unity can reduce difficulty through clever shortcuts or higher-order syntheses. Difficulty also stands in persistent tension with explanatory power, since excessive technical intricacy can obscure the very understanding that explanation demands.     
\item Originality reliably generates fruitfulness, as novel methods tend to be adaptable across domains. The pursuit of unity frequently stimulates originality, and an original framing can enhance explanatory power by dissolving prior obscurities.   
\item Fruitfulness is orthogonal to difficulty—simple proofs can be enormously fruitful—and often depends on originality for its generative reach; it overlaps with unity when cross-domain links open new applications, though it can also remain confined to a single field.    
\item Unity generally demands originality, since recognizing cross-domain connections is an inventive act; it promotes fruitfulness by enabling method-transfer across fields and supports explanatory power by supplying structural overview, yet a proof may unify without explaining (if the connection is technically opaque) or explain without unifying (if it clarifies a single domain through internal insight).
\item Explanatory power often stands in inverse relation to difficulty, as clarity typically requires stripping away unnecessary complexity; it can arise without originality or fruitfulness, and it is frequently delivered through unity, overlapping with both in making a theorem's content intelligible. 
\end{itemize}

The criteria thus form a family of features that cluster in paradigmatic cases, and their value lies not in any single necessary condition but in the structured comparisons they enable across different proofs.
 
\subsection{Are the Criteria Necessary and Sufficient?}\label{sec 3.7}

Consistent with the family-resemblance view introduced in Section~\ref{literature review}, we do not claim that any of the five criteria is strictly necessary or individually sufficient for a proof to be considered deep. A proof that lacks some of difficulty, originality, fruitfulness, unity, or explanatory power might still be regarded as deep, while one that fails to satisfy any of the five is unlikely to be considered deep by reasonable standards. A proof that scores highly on originality and fruitfulness, for instance, may be judged deep even if its difficulty, unity, or explanatory power are only moderate.\footnote{Two further candidates sometimes mentioned in the literature—beauty and influence—deserve brief mention. Beauty, we suggest, is an experiential by-product that arises when several of our criteria converge (e.g., a proof that unifies, explains, and is original is often felt as beautiful), not a conceptually independent dimension of depth (see Inglis and Aberdein [2015]; Thomas [2017, fn. 6]). Influence, while undeniably valuable, largely supervenes on fruitfulness, unity, and originality, and thus risks redundancy as a separate criterion.}
 
Developing a comprehensive list of criteria to measure the depth of mathematical proofs is challenging, if not impossible. Our five-criteria framework establishes a robust heuristic for analysis, acknowledging the distinctions and interconnections among the criteria. While the criteria are neither necessary nor sufficient individually, they offer a balanced and structured vocabulary for comparing proofs and for articulating why some proofs are widely regarded as deep—with Wiles's proof serving as a paradigmatic case, as we will show in Section~\ref{sec4}. 

As a family-resemblance concept, depth may manifest in patterns very different from Wiles's case; our framework is therefore not a general theory of proof depth, but a case-specific analytical grid calibrated to the paradigmatic example that motivates it. For other proofs, different criteria may be more salient, and the framework should be adapted accordingly.

\subsection{The Issue of Objectivity and Subjectivity}\label{sec 3.8}

A recurring theme in the literature is whether judgments of depth can be objective.\footnote{Lange [2014] argues that mathematical depth is context-dependent yet objective, while Urquhart [2014] concludes that depth is relative to audience and historical context.} Three main perspectives emerge: the objectivist view, which holds that depth is an objective feature of proofs revealing real patterns independent of human perception; the subjectivist view, which treats depth as a tool for organizing mathematics on the basis of our interests and cognitive limitations; and the intermediate view, which allows that although depth may depend on our interests and abilities, it can still be objectively assessed using specific criteria (Ernst et al. [2014]). Our framework adopts this intermediate view.  For instance, explanatory power may be judged differently by number theorists (who value high-level unification) and arithmetic geometers (who demand fine-grained intelligibility), yet both judgments are evidence-responsive. By treating explanatory power as graded and multi-faceted, and by allowing that different epistemic values may lead experts to weigh the same evidence differently, the framework provides a structured space for reasoned debate without collapsing into relativism.

Each of the five criteria is context-dependent yet contains objective elements that can be assessed within fixed constraints. Difficulty is influenced by cognitive limitations and historical context—a proof in advanced category theory may daunt a combinatorialist but be routine for a topologist—yet certain aspects, such as formal complexity, admit objective measurement. 
Originality is relative to a community's state of knowledge, but within a fixed historical context one can objectively assess whether the ideas introduced had already been conceived. Fruitfulness depends on our interests and abilities, yet some consequences—such as the number of theorems a proof generates—exist independently of human perception. Unity involves connections between domains that can exist independently of perception, though what counts as ``distinct'' is subject to categorization; given a fixed set of domains, the connections can be objectively assessed. Explanatory power varies across individuals and communities—what one mathematician finds illuminating another may find opaque—and there is no universal standard for a ``good'' explanation, though some aspects (clarity of logic, characterizing properties) admit partial objectification.
Judgments of proof depth, then, are neither wholly objective nor merely subjective; they occupy an intermediate ground.

The depth of a proof is thus evaluated through a combination of criteria, some rooted in mathematics itself and others influenced by human perspectives. Some proofs once considered deep may no longer hold that status today (e.g., the proof of the irrationality of \(\sqrt{2}\) was once profound but is now straightforward). Yet this does not imply that depth is entirely subjective: while academic evaluations contain subjective elements, the mathematical evidence supporting them can be objective. Our intermediate perspective thus transcends simplistic objectivist/subjectivist dichotomies by recognizing the historical and communal nature of judgments on proof depth while still upholding criteria for reasoned debate. 

\section{Assessing the Depth of Wiles's Proof}\label{sec4}

We evaluate the depth of Wiles's proof of FLT using the five criteria introduced above. By ``synthesis'' we mean two things: first, that the proof strongly satisfies each criterion individually; second, and more fundamentally, that the criteria are not independent here. The features that make the proof original also contribute to its unity, and its explanatory power derives partly from the way it unifies disparate fields. Our central thesis is that Wiles's proof exemplifies a non-additive interaction of these dimensions—a mutually reinforcing structure, not a checklist of independent virtues. The following subsections (Sections \ref{sec on difficulty}-\ref{sec on EP}) examine the evidence for each criterion in turn, highlighting the most significant mathematical considerations without attempting exhaustive justification.

\subsection{The Difficulty of Wiles's Proof}\label{sec on difficulty}

Historical evidence, while not itself constitutive of difficulty, strongly indicates the proof's conceptual and technical demands. FLT remained unsolved for 358 years despite sustained efforts by generations of mathematicians, including substantial partial results for many exponents and related developments such as Faltings's Theorem (1983) and the Heath-Brown–Granville ``almost all" result (1984).\footnote{For surveys of these partial results, see Buhler et al. [1993], Faltings [1983] and McLarty [2024].} Yet none of these advances yielded a complete proof or supplied the methods Wiles would later deploy—the critical link to Frey curves had not yet been discovered. The problem was even listed in the Guinness World Records as the ``hardest mathematical problem'' prior to its resolution.
  
The proof's technical difficulty is evident from its reliance on diverse fields—number theory, algebraic geometry, complex analysis, Galois representations, deformation theory, and commutative algebra—each demanding deep specialized expertise. Wiles deployed sophisticated techniques including Galois deformation theory, Hecke algebras, the \(R=T\) strategy, the Taylor–Wiles method, and the ``3–5 switch". The proof spans ninety pages with a 98-entry bibliography, and the discovery of a significant gap during initial peer review further underscores its intricacy.

Conceptually, the proof marks a revolutionary departure. Traditional attempts confined themselves to elementary or algebraic number theory; Wiles instead reframed FLT as a special case of the Taniyama–Shimura–Weil conjecture, establishing a profound structural unity between elliptic curves and modular forms. This paradigm shift proved indispensable—no proof of FLT has been found without it. The methods he introduced to prove the modularity theorem—the modularity lifting theorem, the \(R=T\) strategy, and the Taylor–Wiles method—have since become foundational.

\subsection{The Originality of Wiles's Proof}\label{sec on originality}

Wiles's proof exemplifies originality through both invention and creative synthesis. He transformed FLT from an isolated problem into a corollary of the Taniyama–Shimura–Weil conjecture, linking it to a deeper unifying principle.\footnote{The crucial link between FLT and the conjecture was established by Frey and Ribet, but Wiles was the first to prove the conjecture with the aim of proving FLT.} This reframing is itself a hallmark of originality.

In proving the modularity lifting theorem, Wiles reformulated the modularity problem via the \(R=T\) strategy, turning it into a manageable ring-theoretic question. The Taylor–Wiles method—introducing auxiliary primes, proving isomorphisms in a smoothed setting, and recovering the original via a limiting argument—represents another technical innovation that has become fundamental in number theory.
 
His use of the ``3--5 switch'' to initiate the induction argument further demonstrates originality: rather than proving modularity directly for the reducible mod-3 case, he related it to the mod-5 representation, where modularity was easier to establish, thereby securing the base case.

Beyond these innovations, Wiles's proof achieves a remarkable creative synthesis, integrating decades of work across number theory, algebraic geometry, complex analysis, commutative algebra, and deformation theory. A pivotal instance was his combination of Iwasawa theory and the Kolyvagin–Flach method. Individually, both techniques proved inadequate for bounding the Selmer group—Iwasawa theory failed in the non-abelian context, and the Kolyvagin–Flach Euler system lacked sufficient control. Wiles recognized that their deficiencies were mutually offsetting; by constructing a hybrid framework combining the structural reach of Iwasawa theory with the precise machinery of Kolyvagin–Flach, he leveraged the strengths of each to compensate for the other's weaknesses. This novel integration yielded the critical bound, closed the gap that had stalled the proof for over a year, and completed the modularity theorem for semistable elliptic curves.
 
\subsection{The Fruitfulness of Wiles's Proof}\label{sec on Fruitfulness}

The fruitfulness of Wiles's proof consists in the breadth of its consequences and the applicability of its methods. Unlike FLT itself—easy to state but limited in implications (McLarty [2024])—the Taniyama–Shimura–Weil conjecture has numerous corollaries. 

First, Wiles's proof provides a powerful framework for translating problems about elliptic curves into questions about modular forms and vice versa, a bridge with far-reaching applications that has since been extended well beyond the semistable case.

Second, the methods introduced in Wiles's proof have inspired subsequent generations to tackle more complex problems. His proof established a powerful modular lifting technique, laying the groundwork for a complete proof of the Taniyama–Shimura–Weil conjecture—a development that has proved crucial for further advances in mathematics. The methods Wiles introduced were rapidly generalized: within six years, Conrad–Diamond–Taylor–Breuil extended the modularity theorem to all elliptic curves over the rationals, and subsequent work has extended the scope to real quadratic fields and beyond.

Third, Wiles's proof establishes the $R=T$ strategy, which reframes modularity problems as ring isomorphisms. This strategy translates the abstract question of whether a Galois representation is modular into the concrete algebraic problem of proving that the deformation ring $R$ is isomorphic to the Hecke algebra $T$. The $R=T$ strategy links properties of Hecke algebras to those of Galois representations using commutative algebra. It has since become a powerful tool in algebraic number theory and a standard method for proving modularity lifting theorems [McLarty, 2024].

Fourth, the Taylor-Wiles method introduces auxiliary primes to ``smooth" a deformation problem, proving an isomorphism in a more flexible setting, and then employs a limiting process to return to the original problem. This method has become a fundamental technique in the Langlands program, essential for proving many modularity lifting theorems for higher-dimensional Galois representations.\footnote{Edward Frenkel describes the Langlands program as ``a grand unified theory of mathematics", aiming to construct a comprehensive framework that connects number theory with harmonic analysis and geometry [Gelbart, 1997].}

Fifth, building on earlier work, Wiles developed Galois deformation theory into a central tool. This framework studies families of Galois representations by classifying all possible lifts of a base representation through the ``deformation ring''. This theory is now indispensable in modern number theory, serving as the primary language for formulating and proving results about the modularity of Galois representations, and it is crucial for ongoing work in the Langlands program.
  
In summary, Wiles's proof acts as a catalyst for ongoing mathematical exploration by establishing deep connections between elliptic curves and modular forms. These connections offer a versatile framework that has resolved numerous mathematical problems and opened new research avenues. Moreover, the methods employed—the modularity lifting technique, the $R=T$ strategy, the Taylor-Wiles method, and Galois deformation theory—have provided the mathematical community with a powerful and lasting toolkit.

\subsection{The Unity of Wiles's Proof}\label{sec on Unity}

The unity of Wiles's proof manifests at three interconnected levels: unifying mathematical objects, connecting mathematical fields, and situating the result within the Langlands program.

First,  the proof establishes a connection between elliptic curves and modular forms, encoding the arithmetic data of an elliptic curve—such as the number of solutions modulo various primes—in the Fourier coefficients of its associated modular form. This allows problems in one domain to be translated into the other and solved there. The $R=T$ theorem integrates diverse mathematical objects (e.g., elliptic curve, Galois representation, the deformation ring $R$, and the Hecke algebra $T$) into a cohesive ring-theoretic framework, showing that $R \cong T$ via the Taylor–Wiles method—a powerful tool that simultaneously manages both rings to establish the isomorphism.\footnote{All information about the elliptic curve and its Galois representations is encoded in the universal deformation ring $R$, while the corresponding modular forms are captured by the Hecke Algebra $T$.} 

Second, the proof exemplifies interdisciplinary synthesis, revealing connections among number theory, geometry, analysis, and algebra: 
\begin{itemize}
  \item From arithmetic to geometry: The proof transforms the purely arithmetical problem of Fermat's equation ($a^n + b^n = c^n$) into a geometric problem concerning the properties of a constructed elliptic curve (the Frey curve).\footnote{This shift aligns with the historical association between geometric proofs and depth in number theory, noted by mathematicians such as Riemann [1859], Hilbert and Hurwitz [1890], and Poincar\'{e} [1901].}
\item From geometry to analysis: Wiles's proof of the modularity theorem connects geometry and analysis, demonstrating an intrinsic link between semi-stable elliptic curves and modular forms. 
\item The algebraic machinery: Wiles's proof relies on algebraic structures, particularly Galois representations, which serve as a ``translation device" that converts the abstract symmetries of number fields into the concrete language of linear algebra. This approach forges a link between number-theoretic objects and algebraic structures, utilizing Mazur's theory of Galois representation deformations [McLarty, 2024].
\item Integration of concepts from various fields: Wiles's proof synthesizes concepts from number theory (Fermat's Last Theorem, elliptic curves, Galois groups), algebraic geometry (schemes, cohomology), complex analysis (modular forms), commutative algebra (Hecke algebras, deformation rings), and representation theory (Galois representations).          
\end{itemize}

Third, the proof exemplifies the Langlands program, a broader vision for unifying mathematics.\footnote{Scholars across various fields—such as number theory, harmonic analysis, geometry, representation theory, and mathematical physics—have engaged with the Langlands program, revealing similar unified phenomena even though they focus on different mathematical objects.} 
The Langlands program encompasses profound conjectures about the connections between number theory, analysis, and geometry. The modularity theorem itself is a non-abelian case of more general conjectures proposed by R. Langlands, which seek to associate an automorphic form or representation (a generalization of a modular form) with objects in arithmetic algebraic geometry, such as every elliptic curve over a number field. The proof of the modularity theorem serves as a crucial step toward realizing the broader goals of the Langlands program. By connecting Galois representations generated from elliptic curves with those from modular forms through the \(R = T\) theorem, the proof exemplifies the unity that the Langlands program aims to achieve.

\subsection{The Explanatory Power of Wiles's Proof}\label{sec on EP}

The explanatory power of Wiles's proof lies in two key aspects. 
First, it situates FLT not as an isolated anomaly but as a specific consequence of the deeper modularity theorem. It clarifies the truth of FLT by showing that a counterexample would violate a fundamental principle of mathematics, thereby reframing the question from ``Are there any numbers that satisfy this equation?'' to ``Can a coherent mathematical object like the Frey curve exist?'' The answer is no: a counterexample would create a ``monster'' curve that undermines the essential connection between elliptic curves and modular forms established by the modularity theorem. 

This approach not only resolved FLT but also advanced the understanding of elliptic curves and modular forms. Furthermore, Wiles' proof identifies its characterizing property in the sense of Steiner [1978], revealing that the absence of solutions to Fermat's equation ($a^n + b^n = c^n$) for $n\geq 3$ stems from the fact that a counterexample would correspond to a non-modular semi-stable elliptic curve, contradicting the modularity theorem. Thus, Wiles' proof re-contextualizes FLT—from a mere curiosity about integers to a necessary outcome of structural unity.
 
Second, Wiles' proof provides a generative blueprint for understanding through its core mechanism, the $R=T$ strategy. This strategy elucidates the effectiveness of the modularity lifting technique.  The isomorphism between the deformation ring  $R$  and the Hecke algebra  $T$, established by the $R=T$ theorem, reveals the underlying pattern linking elliptic curves and modular forms. The Taylor-Wiles method further clarifies how to establish this isomorphism between $R$ and $T$.
      
\subsection{The Synthetic Character of Depth}\label{sec 5.7}
 
The preceding subsections show Wiles's proof satisfies each of the five criteria to a high degree. Yet the depth of this proof is not the sum of five independently obtained high scores. The criteria do not merely co-occur; they mutually reinforce one another in a non-additive, structural way. In a merely additive sense, a proof could score highly on each criterion without any internal relation among them. In Wiles's proof, by contrast, the criteria are not independent dimensions.

To make this mutual reinforcement concrete, consider three specific couplings in Wiles' proof:
\begin{itemize}
\item Unity enables explanatory power (Section \ref{sec on Unity}–\ref{sec on EP}): By reducing FLT to the Modularity Theorem, the proof does not merely verify the theorem but explains why a counterexample cannot exist. Without this unifying strategy, no such explanation would be available.
\item Unity and originality together generate fruitfulness (Section \ref{sec on originality}–\ref{sec on Unity}): The $R = T$ strategy is both a unifying framework (connecting deformation rings and Hecke algebras) and an original invention. Precisely because it achieves this structural unification in a novel way, it generalizes far beyond FLT to become a standard tool in the Langlands program. Neither unity alone nor originality alone would have sufficed; their combination produces fruitfulness. 
\item Originality explains the nature of the proof's difficulty (Section \ref{sec on difficulty}–\ref{sec on originality}): The conceptual difficulty of Wiles' proof arises precisely from the need for genuine breakthroughs in ideas and methods. 
\end{itemize}

These couplings are what we mean by `synthesis' – a non-additive interdependence, not merely high scores on independent dimensions.

One might object that if mathematical depth is a family-resemblance concept, then attributing Wiles's depth specifically to his integrative synthesis seems to reintroduce an essence. This objection is misguided. The synthesis in question is not a definitional property of depth but a \emph{paradigm}—a case in which the five hallmark features (difficulty, originality, fruitfulness, unity, and explanatory power) converge and become mutually reinforcing. Other deep proofs may exhibit different convergence patterns. Wiles's proof stands out not because it satisfies a fixed checklist, but because these dimensions cohere into a single, tightly integrated whole. That whole represents one node in the family-resemblance network, not the essence of depth.

In sum, the depth of Wiles' proof resides in the non-additive, mutually reinforcing relations among the five criteria. This synthetic character distinguishes Wiles' achievement from a merely difficult, merely original, or merely unifying proof. It is precisely because the criteria interact and reinforce each other that the proof is widely regarded as a paradigm of mathematical depth.

\section{Conclusion}\label{sec5}
Adopting the family-resemblance perspective, we have treated proof depth not as a definable essence but as a cluster of overlapping, paradigmatically co-occurring features. Our five-criterion framework—Difficulty, Originality, Fruitfulness, Unity, and Explanatory Power—operationalizes this model as a structured heuristic rather than a checklist of necessary and sufficient conditions. This synthesis moves beyond the existing literature by integrating disparate evaluative measures into a coherent analytical grid.
    
The criteria are not merely discrete lenses; their interrelations are constitutive of proof depth itself. While each captures a distinct facet, their overlap is non-redundant: omitting any one dimension risks a caricature of mathematical significance. Depth emerges synergistically—excellence across all five marks a proof as deep, whereas isolated virtuosity in one or two, however brilliant, falls short.

Applied to Wiles's proof, this framework yields a substantive philosophical account. Wiles's achievement exemplifies depth precisely because it satisfies all five criteria in mutual reinforcement, rather than through unilateral virtuosity. We do not claim exhaustiveness, yet the framework effectively demarcates the core evaluative dimensions of mathematical proofs, reconciling their distinctness with their interdependence.
 
The framework serves as a robust tool for paradigmatic case studies—as illustrated by our analysis of Wiles's proof—offering a structured vocabulary for philosophical and methodological reflection. Although qualitatively articulated and eschewing quantification, this is not a methodological deficit but a philosophical necessity: mathematical depth is an inherently evaluative, context-sensitive property that resists metric reduction. Far from hindering consistent application, this very openness enhances the framework's adaptability for probing borderline and contested judgments. Thus, our work does not purport to offer a final taxonomy but rather a fertile starting point for a more nuanced discourse on what constitutes depth in mathematical proofs.

\section*{Acknowledgments}
We are grateful to the editor and the anonymous reviewers for their valuable time and efforts. We particularly appreciate the reviewers' insightful comments and valuable suggestions, which have substantially strengthened this manuscript. 

\section*{References}

Arana, A. [2015]: `On the depth of Szemer\'{e}di's theorem', \emph{Philosophia Mathematica} (3) 23, 163–176. DOI: \url{https://doi.org/10.1093/philmat/nku036} 
 
Arana, A. [2025]: \emph{Elements of Purity}, Cambridge University Press.

Buhler, J.; Crandell, R.; Ernvall, R.; Mets\"{a}nkyl\"{a}, T. [1993]: `Irregular primes and cyclotomic invariants to four million', \emph{Mathematics of Computation}, Vol. 61, No. 203, pp. 151-153. DOI: \url{https://doi.org/10.1090/S0025-5718-1993-1197511-5}

Cheng, Y. [2022]: `On the depth of G\"{o}edel's incompleteness theorems', \emph{Philosophia Mathematica}, Volume 30, Issue 2, pp. 173-199. DOI: \url{https://doi.org/10.1093/philmat/nkab034}

Cornell, G.;  Silverman, J. H.  and Stevens, G.  [1997]: \emph{Modular forms and Fermat's Last Theorem}, Springer-Verlag.

Darmon, H.; Diamond, F.; and Taylor, R. [1997]: `Fermat's Last Theorem', In Coates, J. and Yau, S.-T., editors, \emph{Elliptic Curves, Modular Forms and Fermat's Last Theorem}, pages 2–140. International Press, Somerville Massachussetts.

Demyanenko, V. [1971]: `The points of finite order of elliptic curves', \emph{Acta Arith}, 19:185–194
 
Ernst M.; Heis, J.; Maddy, P.;  McNulty, B.M.; Weatherall, O. J. [2014]: `Foreword to Special Issue on Mathematical Depth', `Afterword to Special Issue on Mathematical Depth', \emph{Philosophia Mathematica}, Volume 23, Issue 2. DOI: \url{https://doi.org/10.1093/philmat/nkv003}; DOI: \url{https://doi.org/10.1093/philmat/nkv002}

Faltings G. [1983]: `Endlichkeitssatze f\"{u}r abelsche Variet\"{a}ten\"{u}ber Zahlk\"{o}rpern', \emph{Invent Math} 73:349–366. DOI: \url{https://doi.org/10.1007/BF01388432}

Frey, G. [1986]: \emph{`Links between stable elliptic curves and certain Diophantine equations'},
Number 1 in Annales Universitatis Saraviensis. Series Mathematicae. Saarlande University, Saarbr\"{u}cken. 

Gelbart, S. [1997]: `Three lectures on the modularity of $\overline{\rho}_{e,3}$ and the Langlands reciprocity conjecture', 
In: Cornell G, Silverman J, Stevens G (eds), \emph{Modular forms and Fermat's last theorem}, Springer,
New York, pp 155–191.  

Gray, J. [2014]: `Depth—A Gaussian Tradition in Mathematics', \emph{Philosophia Mathematica},  Volume 23, Issue 2, 177–195. DOI: \url{https://doi.org/10.1093/philmat/nku035}

Hafner, J., and Mancosu, P. [2005]: `The varieties of mathematical explanation', in P. Mancosu (ed.), Visualization, Explanation and Reasoning Styles in Mathematics, Springer, pp. 215–250. 

Hellegouarch, Y. [1971]: `Points d'ordre fini sur les courbes elliptiques', \emph{CRAS Paris},
273:54–43.

Hilbert, D. and Hurwitz, A. [1890]: `\"{U}ber die diophantischen Gleichungen vom
Geschlecht Null', \emph{Acta Mathematica}, 14:217–224. DOI: \url{https://doi.org/10.1007/BF02413323}
 
Holden, H.; Piene, R. [2014]: \emph{The Abel Prize 2008–2012}, Springer Science+Business Media B.V.
 
Inglis, M.; Aberdein, A. [2015]: `Beauty Is Not Simplicity: An Analysis of
Mathematicians' Proof Appraisals', \emph{Philosophia Mathematica}, Volume 23, Issue 1, pp. 87–109. DOI: \url{https://doi.org/10.1093/philmat/nku014}

Lange, M. [2014]: `Depth and Explanation in Mathematics', \emph{Philosophia Mathematica}, Volume 23, Issue 2, 1–19. DOI: \url{https://doi.org/10.1093/philmat/nku022} 

Mazur, B. [1977]: `Modular curves and the Eisenstein ideal', \emph{Publ. Math. IHES},
47:133–86. DOI: \url{https://doi.org/10.1007/BF02684339}
 
Mazur, B.; Goldfeld, D. [1978]: `Rational isogenies of prime degree', \emph{Inventiones mathematicae},
44:129–62. DOI: \url{https://doi.org/10.1007/BF01390348}

Mazur, B. [1991]: `Number theory as gadfly', \emph{American Mathematical Monthly},
98: 593–610. DOI: \url{https://doi.org/10.1080/00029890.1991.11995762}
 
McLarty, C. [2024]: `Fermat's last theorem', In Sriraman, B., editor, \emph{Handbook of
the History and Philosophy of Mathematical Practice}, Springer, pp. 2011–2033.

Poincar\'{e}, H. [1901]: `Sur les propri\'{e}t\'{e}s arithm\'{e}tiques des courbes alg\'{e}briques',  \emph{Journal des Math\'{e}matiques}, 7:161–233. 

Ribet, K. [1990]: `On modular representations of $Gal(Q/Q)$ arising from modular forms',  \emph{Inventiones Mathematicae}, 100 (2): 431–476. DOI: \url{https://doi.org/10.1007/BF01231195} 

Riemann, B. [1859]: Ueber die Anzahl der Primzahlen unter einer gegebenen
Gr\"{o}sse, Handwritten note, reprinted in Harold Edwards Riemann's zeta function,
Dover Publications 2001, among other places.

%Singh, S. [1997]: \emph{`Fermat's Last Theorem'}, Fourth Estate Ltd, ISBN 1-85702-521-0.
 
Stillwell, J. [2014]: `What Does `Depth' Mean in Mathematics?', \emph{Philosophia Mathematica}, Volume 23, Issue 2, 215–232. DOI: \url{https://doi.org/10.1093/philmat/nku033}

Steiner, M. [1978]: `Mathematical Explanation', \emph{Philosophical Studies}, 34(2), 135–151. DOI: \url{https://doi.org/10.1007/BF00354494}

Tappenden, J. [2012]: `Fruitfulness as a Theme in the Philosophy of Mathematics', \emph{The Journal of Philosophy}, Vol. 109, No. 1/2, pp. 204-219. DOI: \url{https://doi.org/10.5840/jphil20121091/27}
 
Thomas, R. S. D. [2017]: `Beauty is not all there is to Aesthetics in Mathematics', \emph{Philosophia Mathematica}, Volume 25, Issue 1, pp. 116–127. DOI: \url{https://doi.org/10.1093/philmat/nkw019}
 
Urquhart, A. [2014]: `Mathematical Depth', \emph{Philosophia Mathematica},  Volume 23, Issue 2, 233–241.  DOI: \url{https://doi.org/10.1093/philmat/nkv004}
 
Wiles, A. [1995]: `Modular elliptic curves and Fermat's Last Theorem', \emph{Annals of Mathematics}, 141:443–551. DOI: \url{https://doi.org/10.2307/2118559}
 
Wittgenstein, L. [2009]: Philosophical Investigations (G. E. M. Anscombe, P. M. S. Hacker, \& J. Schulte, Trans., 4th ed.). Wiley-Blackwell.
     
\end{document}